\documentclass[a4paper]{amsart}
\usepackage[latin1]{inputenc}
\usepackage{amsmath}
\usepackage{amsfonts}
\usepackage{amssymb}
\usepackage{enumerate}
\usepackage{mathrsfs}
\usepackage[all]{xy}
\usepackage{epic}

\newcommand\bk{{\mathbb K}}
\newcommand\bp{{\mathbb P}}
\newcommand\bff{{\mathbb F}}
\newcommand\br{{\mathbb R}}
\newcommand\bc{{\mathbb C}}
\newcommand\bq{{\mathbb Q}}

\newcommand\bz{{\mathbb Z}}

\newcommand\ca{{\mathcal A}}
\newcommand\calp{{\mathcal P}}

\newcommand\scl{{\mathscr L}}

\newcommand\scp{{\mathscr P}}

\newcommand\codim{{\text{codim}}}

\DeclareMathOperator\corank{corank}

\DeclareMathOperator\Aut{Aut}

\newtheorem{thm}{Theorem}[section]
\newtheorem{prop}[thm]{Proposition}
\newtheorem{cor}[thm]{Corollary}
\newtheorem{lema}[thm]{Lemma}

\theoremstyle{remark}
\newtheorem{rmrk}[thm]{Remark}

\theoremstyle{definition}
\newtheorem{dfn}[thm]{Definition}
\newtheorem{ntc}[thm]{Notation}
\newtheorem{ejm}[thm]{Example}

\numberwithin{equation}{section}

\title{Resonance varieties, admissible line combinatorics and combinatorial pencils}
\author{Miguel \'Angel Marco Buzun\'ariz}
\address{Departamento de matem\'aticas \\
 Universidad de Zaragoza\\
 C/ Pedro Cerbuna 12\\
 50009 Zaragoza\\
 Spain}
 \email{mmarco@unizar.es}

\begin{document}

\bibliographystyle{amsplain}

\begin{abstract}
In this paper we define the combinatorial analogous of a pencil, and show its relationship with the concept of admissibility. Such an object is usefull to study the isomorphisms between fundamental groups of the complements of line arrangement with the same combinatorial type. This definition generalizes the idea of net given by Yuzvinsky and others. The main theorem in this paper states that there is a correspondence between components of the resonance variety and combinatorial pencils.
\end{abstract}

\maketitle
\section*{Introduction}
The problem of the relationship between the topology of a line arrangement in $\bc \bp ^2$ and its combinatorial structure has been one of the most studied in the theory of hyperplane arrangements. After the work of Arnol'd and Brieskorn, Orlik and Solomon showed that the cohomology algebra of the complement is determined by the intersection lattice. In particular this implies that the first homology group depends only on the number of lines. Rybnikov exhibited in \cite{rybnikov} the existence of two combinatorially equivalent line arrangements $\scl_1,\scl_2$, whose complements have non-isomorphic fundamental groups. His approach had essentially two parts: On one hand, he showed that an isomorphism between the fundamental groups should send meridians to meridians in homology. On the other hand, he could distinguish both arrangements using an invariant under such isomorphisms. In the first part is where most details were missing.

In \cite{nuestro-ryb} a detailed proof of Rybnikov's result s given. In this paper we generalize the concepts and methods introduced in \cite{nuestro-ryb}. The concept of admissible classes is introduced. It is strongly related to other concepts spread in the literature, such as the idea of component of the resonance variety, characteristic varieties, neighborly partitions and nets (see \cite{falk-arrco,libgober-yuzvinsky,yuzvinsky-nets}). In particular, each maximal admissible class corresponds to a component of the resonance variety. The main result (Theorem~\ref{admhaz}) of this paper shows that each admissible class determines a combinatorial pencil. This equivalence has not only theoretical importance, but also computational: checking the admissibility of a combinatorics involves solving a system of quadratic equations, whereas doing the same for a combinatorial pencil only involves solving linear systems. 

Section \ref{definiciones} contains the definitions of admissibility for maps, classes and combinatorics. Some examples of admissible combinatorics are shown. It also contains the definition of combinatorial pencil, which  generalizes the idea of net, and is a particular case of neighborly partitions. The proof of the main theorem is provided in Section \ref{decfibras}. It uses the Vinberg classification of real matrices, which is included for completeness. The relationship between admissible classes and components of the resonance variety is also shown in this section.
Section \ref{permutaciones} describes a method to study the isomorphisms between fundamental groups of combinatorially equivalent  line arrangements, using the previous concepts. This method involves studying  the permutation induced in the set of maximal admissible classes, $Adm(\scl,\scp)$. Such a permutation must preserve some structure in $Adm(\scl,\scp)$, such as the triangles of admissible classes. The concept of triangle of admissible classes is defined and shown to be invariant under permutation induced by isomorphisms of the fundamental group. This method is used in Section \ref{ejemplo} to show the homological rigidity of a certain ten line combinatorics that admits real realizations.
%
Section \ref{dualidad} includes, for the sake of completeness, a proof of the well-known duality between the second level of the lower central series of the fundamental group of a line arrangement, and its Orlik-Solomon algebra. This fact is important because it opens the door to use this approach in the attempt to study the completeness of the Orlik-Solomon algebra as invariant of the combinatorial type. The aforementioned method can be used to prove that certain combinatorics are homologically rigid, which is a sufficient condition to conclude that there is no other combinatorial type with isomorphic OS algebra.

\section{Preliminary definitions}
\label{definiciones}

\begin{dfn} A \textbf{line combinatorics} is a finite set $\scl:=\{l_1,\ldots , l_n\}$ together with a subset $\scp\subseteq\calp(\scl)$ satisfying the following properties:
\begin{enumerate}[i)]
\item $\# p>1, \ \ \forall p\in\scp$.
\item for every $l_i,l_j\in\scl$, $l_i\neq l_j$, there exists a unique $p\in\scp$ such that $l_i, l_j\in p$. This element will be called the intersection of $l_i$ and $l_j$, and will be denoted by $l_i \cap l_j$.
\end{enumerate}

The elements of $\scl$ and $\scp$ will be called \textbf{lines} and \textbf{points} respectively.
\end{dfn}
This definition captures the incidence properties of the set of lines and the set of points (identifying a point with the set of lines that pass through it) of a line arrangement in a projective plane. Given a line arrangement in the complex projective plane, some of the invariants of its topology depend only on the combinatorics; so we will reffer to them as invariants of the combinatorics.
In the following we will assume that we have fixed a line combinatorics of $n$~lines. In this context, we will define $H$ as the quotient of the lattice $\bz ^n$ by the sublattice generated by the vector $(1,\ldots,1)$. We fix an order in $\scl$ that allows us to establish a bijection between the lines $\{l_1,\ldots,l_n\}$ and the elements of the canonical generating system $\{e_1,\ldots, e_n\}$ of $H$, (these are the classes in $H$ of the canonical basis of $\bz ^n$). For the sake of simplicity, $e_l$ with $l\in\scl$ will also denote a canonical generator of $H$.  If there existed a realization of the combinatorics in the complex projective plane, $H$ would be canonically isomorphic to the first homology group of the complement.
\begin{dfn}
Let $k$ be a positive integer greater than $2$. A \textbf{k-admissible map} is a $\bz$-epimorphism $\alpha :H\rightarrow \bz ^{k-1}$ such that, for every $p\in\scp$ and for every $l_i\in p$, the vectors $\{\alpha(e_i), \sum_{l_j\in p} \alpha(e_j)\}$ are linearly dependent.

Given an admissible map, we define its associated subcombinatorics as the combinatorics whose  set of lines is $\scl_{\alpha}:=\{l_i \mid \alpha(e_i)\neq 0\}$, and its set of points is $\scp_{\alpha}:=\{p \cap \scl_{\alpha} \mid p\in \scp, \#(p \cap \scl_{\alpha}) > 1\}$. That reflects the intuitive idea of ``deleting" the lines $l_i$ where $\alpha (e_i)$ vanishes.

A line combinatorics is said to be \textbf{k-admissible} if it admits an admissible map $\alpha$ such that $\alpha (e_l)\neq 0\ \forall l\in\scl$.

The group $Aut_{\bz}(\bz^{k-1})$ acts on the set of $k$-admissible maps by composition. The  orbits of this action will be called \textbf{k-admissible classes}.

We will say that a $k$-admissible map is maximal if it can not be obtained by composition of a $(k+1)$-admissible map and a $\bz$-epimorphism $\bz^k\rightarrow \bz^{k-1}$. Analogously we will talk about maximal $k$-admissible combinatorics and classes. The set of maximal $k$-admissible classes will be denoted by $Adm_{k}(\scl,\scp)$; and the set of all maximal admissible classes (that is, $\bigcup_{k>2}Adm_{k}(\scl,\scp)$) will be denoted by $Adm(\scl,\scp)$.
\end{dfn}
 
The combinatorics $(\scl_\alpha,\scp_\alpha)$ associated with an admissible map $\alpha$ is invariant of its admissible class.

\begin{ejm}
\label{punto}
Let $l_{1},\ldots ,l_{k}$ be $k$ concurrent lines (with $k\geq 3$). One can define a $k$-admissible map as follows: consider $\{v_1,\ldots ,v_{k-1}\}$ a basis of $\bz^{k-1}$. Now let $\alpha$ be the map given by $\alpha (e_j)=v_j$ for $j=1,\ldots, k-1$, and $\alpha(e_{k})=-v_1-\cdots -v_{k-1}$. This is an admissible map, in fact this $k$-admissible combinatorics is maximal, we will call it of \emph{point type}.
An admissible map whose associated subcombinatorics is of point type is also called an admissible map of \emph{point type}.
\end{ejm}
\begin{ejm}
\label{ceva}
Let $l_{1},\ldots ,l_{6}$ be six lines whose non-double points are $\{l_{1},l_{2},l_{3}\}$, $\{l_{1},l_{5},l_{6}\}$, $\{l_{2},l_{4},l_{6}\}$, $\{l_{3},l_{4},l_{5}\}$. Fixing $\{v_1,v_2\}$ to be a basis of $\bz^2 $, we can define the following 3-admissible map $\alpha$ given by $\alpha(e_{1})=\alpha(e_{4})=v_1$, $\alpha(e_{2})=\alpha(e_{5})=v_2$, $\alpha(e_{3})=\alpha(e_{6})=-v_1-v_2$. It is again easy to check that $\alpha$ is admissible.
The admissible classes of this form will be called of \emph{Ceva type}. As in the previous example, we can talk about subcombinatorics of \emph{Ceva type}.
\end{ejm}
\begin{ejm}
\label{cuerpofinito}
Consider a finite field $\bff$ of $k$ elements, and consider $\bff^2$ the affine plane over $\bff$. In this plane there are $k+1$ directions, and for every direction there are $k$ lines. For each point of the plane passes exactly one line of each direction. So we can construct $(k+1)$-admissible combinatorics as follows: for each direction $D$, construct a combinatorics $(\scl_D,\scp_D)$ using the lines whose direction is $D$, and consider the combinatorics $(\scl,\scp)$ such that $\scl$ is the set of all lines in $\bff^2$, and $\scp$ is the union of all the $\scp_D$ plus the points of $\bff^2$. Then choose a basis $(e_1,\ldots ,e_k)$ of $\bz^k$, and order the directions of the plane. Now map the lines in the $i$'th direction to $e_i$, and the ones in the last direction to $-e_1-\cdots -e_k$. This map will be admissible regardless of the election of $\scp_D$.

In the case of $\bz / 2\bz$, there is only one way of choosing the intersections inside each direction (two lines only have one way to intersect), and the result is the Ceva combinatorics.

In the case of $\bz / 3\bz$, and choosing the lines in each direction to be in general position, the result is the combinatorics of the twelve lines joining the nine flexes of a smooth cubic. This one is called the Hesse combinatorics.

There is also a ``degenerated" Hesse combinatorics, which is not realizable in the complex plane, where all the parallel lines intersect at one point.
\end{ejm}

\begin{ejm}
\label{cevageneralizado}
From the Ceva combinatorics, we can add lines joining every pair of double points. The result is a combinatorics of nine lines, as shown in Figure~\ref{cevas}, and is called \emph{generalized Ceva} combinatorics. In Figure~\ref{cevas} we can see that it is admissible. Note that in this case, there are non-equal proportional vectors.

The realization in the complex (or real) plane is the union of the three special fibers of the pencil of rational nodal quartics, generated by $x^2(y^2+z^2)$ and $y^2(x^2+z^2)$ (up to projective transformation).
\end{ejm}
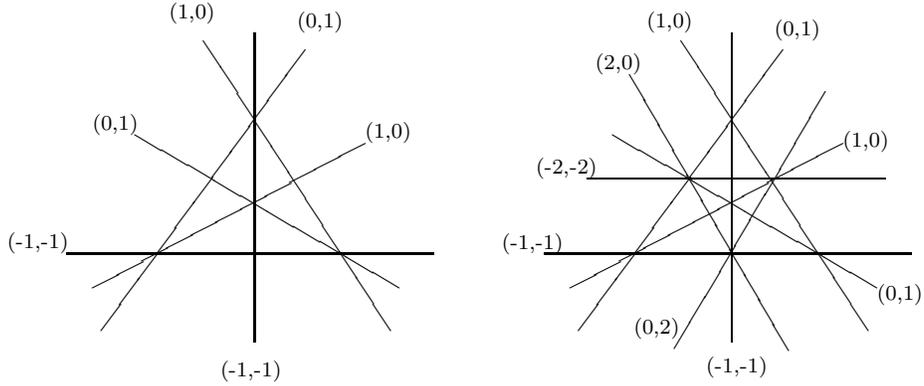
\begin{figure}[h]
\caption{\label{cevas} Ceva and generalized Ceva combinatorics}
\setlength{\unitlength}{0.0006in}
\begingroup\makeatletter\ifx\SetFigFont\undefined%
\gdef\SetFigFont#1#2#3#4#5{%
  \reset@font\fontsize{#1}{#2pt}%
  \fontfamily{#3}\fontseries{#4}\fontshape{#5}%
  \selectfont}%
\fi\endgroup%
{\renewcommand{\dashlinestretch}{30}
\begin{picture}(8044,3448)(0,200)
\drawline(2625,3054)(825,579)
\drawline(1725,3129)(3375,579)
\drawline(2175,3204)(2175,479)
\drawline(525,1254)(3750,1254)
\drawline(750,954)(3150,2229)
\drawline(3450,954)(1125,2304)
\drawline(6825,3054)(5025,579)
\drawline(5925,3129)(7575,579)
\drawline(6375,3204)(6375,479)
\drawline(4725,1254)(7950,1254)
\drawline(4950,954)(7350,2229)
\drawline(7650,954)(5325,2304)
\drawline(7200,2679)(5870,421)
\drawline(6875,404)(5475,2829)
\drawline(5100,1915)(7725,1915)
\put(0,1254){\makebox(0,0)[lb]{{\SetFigFont{8}{10}{\rmdefault}{\mddefault}{\updefault}(-1,-1)}}}
\put(1875,129){\makebox(0,0)[lb]{{\SetFigFont{8}{10}{\rmdefault}{\mddefault}{\updefault}(-1,-1)}}}
\put(1425,3279){\makebox(0,0)[lb]{{\SetFigFont{8}{10}{\rmdefault}{\mddefault}{\updefault}(1,0)}}}
\put(3150,2229){\makebox(0,0)[lb]{{\SetFigFont{8}{10}{\rmdefault}{\mddefault}{\updefault}(1,0)}}}
\put(2550,3204){\makebox(0,0)[lb]{{\SetFigFont{8}{10}{\rmdefault}{\mddefault}{\updefault}(0,1)}}}
\put(750,2304){\makebox(0,0)[lb]{{\SetFigFont{8}{10}{\rmdefault}{\mddefault}{\updefault}(0,1)}}}
\put(6750,3129){\makebox(0,0)[lb]{{\SetFigFont{8}{10}{\rmdefault}{\mddefault}{\updefault}(0,1)}}}
\put(7650,804){\makebox(0,0)[lb]{{\SetFigFont{8}{10}{\rmdefault}{\mddefault}{\updefault}(0,1)}}}
\put(5625,3204){\makebox(0,0)[lb]{{\SetFigFont{8}{10}{\rmdefault}{\mddefault}{\updefault}(1,0)}}}
\put(7350,2154){\makebox(0,0)[lb]{{\SetFigFont{8}{10}{\rmdefault}{\mddefault}{\updefault}(1,0)}}}
\put(6150,160){\makebox(0,0)[lb]{{\SetFigFont{8}{10}{\rmdefault}{\mddefault}{\updefault}(-1,-1)}}}
\put(4350,1250){\makebox(0,0)[lb]{{\SetFigFont{8}{10}{\rmdefault}{\mddefault}{\updefault}(-1,-1)}}}
\put(5510,504){\makebox(0,0)[lb]{{\SetFigFont{8}{10}{\rmdefault}{\mddefault}{\updefault}(0,2)}}}
\put(5175,2829){\makebox(0,0)[lb]{{\SetFigFont{8}{10}{\rmdefault}{\mddefault}{\updefault}(2,0)}}}
\put(4650,1910){\makebox(0,0)[lb]{{\SetFigFont{8}{10}{\rmdefault}{\mddefault}{\updefault}(-2,-2)}}}
\end{picture}
}
\end{figure}

The concept of admissible map can be presented in a way that is more independent of the election of a basis in $H$. Consider the $\bz-$submodule $R\leqslant H\wedge H$ generated by the family $\{e_j\wedge \sum_{i \in p} e_i \mid p\in \scp\ , j\in p\}$ (where $\wedge$ denotes the exterior product).

\begin{prop}
\label{carac-admisible}
An epimorphism $\alpha :H\rightarrow \bz ^{k-1}$ is k-admissible if and only if $\alpha \wedge \alpha (x)=0\ \forall x\in R$.
\end{prop}
\begin{proof}
Let $p$ be a point, $l_j\in p$ and consider a generator $e_j\wedge \sum_{i \in p} e_i$ of $R$. Its image under $\alpha\wedge \alpha$ is $\alpha(e_i)\wedge \sum_{i \in p} \alpha(e_i)$, which is zero if and only if $\alpha(e_i)$ and $\sum_{i \in p} \alpha(e_i)$ are linearly dependent in~$\bz^{k-1}$. 
\end{proof}

The definition of $R$ is motivated by the complement of a realization in the complex projective plane: the second term of the lower central series of its fundamental group is isomorphic to $(H\wedge H) / R$, see~\cite{rybnikov}. This group also appears in the study of the truncated Alexander invariant of the complement of a line arrangement. 
In \cite{nuestro-ryb} this group was used to study the set of isomorphisms of fundamental groups  of both McLane's and Rybnikov's arrangements.

A straightforward consequence of \cite[Prop. 7.2]{libgober-yuzvinsky} is that a line arrangement whose combinatorics is $k$-admissible is a union of fibers of a pencil. In particular, Example~\ref{punto} is trivially the union of some fibers of a pencil of lines; Example~\ref{ceva} is the union of the three singular fibers of a pencil of conics in general position; and Example ~\ref{cevageneralizado}, as explained before, is the union of three non-reduced fibers of a pencil of quartics.  In the following we will define a generalization of the concept of pencil in purely combinatorial terms:
\begin{dfn}
\label{combinatorialpencil}
A \textbf{combinatorial pencil} is a line combinatorics together with a partition $F_1,F_2,\ldots ,F_k$ of $\scl$ and a weight map $w:\scl \rightarrow \bz^+$ such that at any point $p\in\scp$ only one of the following two possibilities occurs:
\begin{enumerate}[i)]
\item $p \subseteq F_i$ for some $i \in \{1,2,\ldots , k\}$.
\item $\forall i \in \{1,2,\ldots ,k\}$, $p\cap F_i\neq\emptyset$ and $\sum_{l\in p\cap F_i} w(l)=\sum_{l\in p\cap F_j} w(l)$ for all $F_i,F_j$.
\end{enumerate}
The points satisfying property ii) will be called the \textbf{base points} of the combinatorial pencil, and the elements $F_i$ of the partition will be called \textbf{fibers}.
\end{dfn}
\begin{rmrk}
\label{penciladm}
From a combinatorial pencil of $k$ fibers $F_1,\ldots , F_k$, we can construct a $k$-admissible map as follows: let $\{v_1,\ldots ,v_{k-1}\}$ be a basis of $\bz^{k-1}$, and $v_k := -v_1- \cdots -v_{k-1}$, now consider the map $\alpha: H \rightarrow \bz^{k-1}$ given by $\alpha (e_i):=w(l_i)\cdot v_j$, where $l_i\in F_j$. The properties of Definition~\ref{combinatorialpencil} imply that $\alpha$ is $k$-admissible. Thus, each combinatorial pencil determines an admissible map; the goal of Section~\ref{decfibras} is to prove its converse. 
\end{rmrk}
\begin{rmrk}
By the aforementioned result \cite[Prop. 7.2]{libgober-yuzvinsky}, in case the combinatorics is realizable, any combinatorial pencil is in fact a geometrical pencil.
\end{rmrk}

\section{Decomposition in fibers}
\label{decfibras}

The goal of this section is to prove that each admissible map determines a combinatorial pencil. In order to do so, we will use some ideas from \cite{libgober-yuzvinsky} and the Vinberg classification of matrices (see \cite{vinberg}), which we will include here for completeness. In particular, we will use \cite[Thm. 4.3]{vinberg}:
\begin{ntc}
Given a vector $u\in \br^n$, we will write $u\geq 0$ (resp. $>$,$\leq$ or $<$) to denote that all its entries are nonnegative (resp. positive, nonpositive or negative).
\end{ntc}
\begin{thm}
\label{vinberg}
Let $A=(a_{i,j})$ be a real $n\times n$ matrix such that:
\begin{itemize}
\item $A$ is indecomposable.
\item $a_{i,j}\leq 0$ for $i\neq j$.
\item $a_{i,j}=0$ implies $a_{j,i}=0$. 
\end{itemize}
Then only one of the following three possibilities hold for both $A$ and its transposed:
\begin{itemize}
\item[(Fin)] $det(A)\neq 0$; there exists $u>0$ such that $Au>0$; $Av \geq 0$ implies $v=0$ or $v>0$.
\item[(Aff)] $\corank (A)=1$; there exists $u>0$ such that $Au=0$; $Av\geq 0$ implies $Av=0$.
\item[(Ind)] there exists $u>0$ such that $Au<0$; $Av\geq 0$, $v\geq 0$ imply $v=0$.
\end{itemize}
\end{thm}

In all this section we will assume that $(\scl,\scp)$ is a maximal $k$-admissible line combinatorics with admissible map $\alpha$. First consider $\chi_{\alpha}:=\{p\in\scp \mid \sum_{l_i\in p} \alpha(e_i)=0\}$. 

If $\#\chi_\alpha=1$ then $(\scl,\scp)$ is of point type, otherwise there exists a line $l_i$ that does not go through the only point in $\chi_\alpha$. Then for every line $l_j$, the corresponding vector $\alpha (e_j)$ is proportional to $\alpha(e_i)$, since both are proportional to $\sum_{l_k\in l_i\cap l_j}\alpha(e_k)$, which is not zero. This contradicts the admissibility of $\alpha$. From Example~\ref{punto} a maximal $k$-admissible combinatorics of point type defines a combinatorial pencil with $k$ fibers (one per line).

From now, we will assume that $\#\chi_\alpha\geq 2$. Consider a graph whose vertices are the lines and whose edges join every two lines that intersect outside $\chi_\alpha$. We have a partition of $\scl$ given by the connected components of this graph; let's denote such a partition by $\Pi$. We can assume that the lines are ordered in a way compatible with $\Pi$ (that is: if $l_i$ and $l_j$ are in the same component, and $i<k<j$, then $l_k$ is also in the component of $l_i$ and $l_j$). 
\begin{rmrk}
\label{matrizq}
This same decomposition is done in \cite{libgober-yuzvinsky} with a slightly different approach. Consider $Q$ the $(n\times n)$ matrix whose entries are $Q_{i,i}:=\#\{p\in \chi_{\alpha} \mid l_i \in p\}-1$ in the diagonal; $Q_{i,j}:=0$ if the intersection of $l_i$ and $l_j$ is in $\chi_{\alpha}$; and $Q_{i,j}:=-1$ otherwise. This is a symmetric matrix that can be decomposed in a direct sum of indecomposable matrices $\oplus_{F\in\Pi}M_F$. This decomposition corresponds to the connected components of the previous graph. It is straightforward to check that $Q_{i,i}\geq 1$ for every $i$.
Another way to define $Q$ is $Q:=J^TJ-U$, where $J$ is the incidence matrix between $\chi_{\alpha}$ and $\scl$, and $U$ is the $\#\scl\times\#\scl$ matrix whose every entry is $1$. 
\end{rmrk}

Note that if two lines $l_i,l_j$ are in the same component of $\Pi$, the vectors $\alpha(e_i)$ and $\alpha(e_j)$ are linearly dependent, and hence if $F\in \Pi$, there exists a primitive vector $v_F\in\bz^{k-1}$ such that $\forall l \in F$, $\alpha(e_l)=w_lv_F$ for some $w_l\in \bz$.

\begin{lema}
Let $F\in \Pi$, all the entries of the weight vector $(w_l)_{l\in F}$ have the same sign. In particular $(w_l)_{l\in F}$ can be chosen to be positive.
\end{lema}
\begin{proof}
Fix a line $l_i$, the following equations hold
\begin{equation}
\label{sistematodos}
\left\{\left.\sum_{l_j\in p}\alpha(e_j)=0\ \right\vert\ p\in \chi_\alpha ,\ l_i \in p\right\}
.
\end{equation}
The following properties hold for the system~\eqref{sistematodos} 
\begin{itemize}
\item $\alpha(e_i)$ appears $Q_{i,i}+1$ times,
\item if $l_i\cap l_j\in \chi_\alpha$, $\alpha(l_j)$ appears exactly once,
\item if $l_i\cap l_j\notin \chi_\alpha$, $\alpha(l_j)$ does not appear.
\end{itemize}
Hence if we substract every equation in~\eqref{sistematodos} from the equation
\begin{equation}
\label{sumacero}
\sum_{j=1}^n\alpha(e_j)=0
,
\end{equation}
we obtain that 
\begin{equation}
\label{anulafilas}
-Q_{i,i}\alpha(e_i)+\sum_{l_i\cap l_j \notin \chi_\alpha}\alpha(e_j)=0
.
\end{equation}
Since for all such $l_j$, $\alpha(e_j)=w_j v$ for a certain $v\neq 0$, the equation~\eqref{anulafilas} can be expressed as
\begin{equation}
\label{anulafilaspesos}
Q_{i,i}w_i-\sum_{l_i\cap l_j \notin \chi_\alpha}w_j=0
,
\end{equation}
which means in particular that the weight vector $(w_{1},\ldots,w_{n})$ is in the kernel of $Q$. Moreover, for every $F\in \Pi$, the weight vector associated with $F$, $(w_l)_{l\in F}$ is in the kernel of the corresponding indecomposable matrix $M_F$. Since all these matrices satisfy the hypothesis of Theorem~\ref{vinberg}, $M_F$ is of one of the three types (Aff), (Fin) or (Ind). We have found a member of its kernel, so it cannot be of (Fin) type. 
Now suppose that, for a certain $G\in \Pi$, the matrix $M_G$ is of (Ind) type. There exists a positive vector $u_G>0$ such that $M_Gu_G<0$. Now, for every $F\in \Pi\setminus \{G\}$, there exist a vector $u_F<0$ such that $M_Fu_F>0$ (if $M_F$ is of (Ind) type), or $M_Fu_F=0$ (if $M_F$ is of (Aff) type). By multiplying each $u_F$ by an adequate positive constant, we can assume that $\sum_{l\in \scl}u_l=0$. Consider the vector $u=(u_l)_{l\in\scl}$ obtained by concatenation of all the $u_F$. Since the sum of the entries of $u$ is zero, $Uu=0$. Then, denoting by $(\bullet,\bullet)$ the standard escalar product, we have:
\begin{equation}
0\leq (Ju,Ju)=(Qu,u)+(Uu,u)=(M_Gu_G,u_G)+\sum_{F\in\Pi\setminus\{G\}}(M_Fu_F,u_F)\leq (M_Gu_G,u_G)<0
\end{equation}
which is a contradiction. We conclude that all the $M_F$ are of (Aff) type, and since the vectors $(w_l)_{l\in F}$ generate the kernel of $M_F$, all its entries must have the same sign.
 
\end{proof}

In particular, since $\alpha(e_1),\ldots,\alpha(e_n)$ generate $\bz^{k-1}$ and $\sum_{i=1}^n e_i=0$, the previous $\{v_F\}_{F\in\Pi}$ is a linearly dependent generating system in $\bz^{k-1}$. So we can conclude that $\#\Pi\geq k$.

We now will recall the definition of the Orlik-Solomon algebra of a line combinatorics (see \cite{orlik-terao}) in order to use some of its properties.
\begin{dfn}
Let $(\scl,\scp)$ be a line combinatorics. $\scl=\{l_1,\ldots ,l_n\}$.
Consider the ($n-1$)-dimensional vector space $E_1$ over a field $\bk$ generated by $x_1,\ldots,x_n$ satisfying the relation $x_1+\cdots +x_n=0$. Let $E$ be the graded exterior algebra of $E_1$ (note that $x_1,\ldots,x_n$ correspond to the generators $e_1,\ldots,e_n$ of $H$ but it is more convenient to use a different notation to distinguish both objects). Now consider the differential $\delta:E_p\rightarrow E_{p-1}$ given by
$$\delta(x_{i_1} \wedge \cdots \wedge x_{i_p})=\sum_{j=1}^p(-1)^{j-1}(x_{i_1}\wedge \cdots \wedge \hat{x}_{i_j}\wedge \cdots \wedge x_{i_p})
$$
The \textbf{Orlik-Solomon algebra} over $\bk$ of $(\scl,\scp)$ is defined as the quotient $\ca$ of $E$ by the ideal generated by $\{\delta(x_{i_1}\wedge \cdots \wedge x_{i_n})\mid n>3\}$, $\{\delta(x_{i_1}\wedge x_{i_2} \wedge x_{i_3})\mid l_{i_1}\cap l_{i_2}=l_{i_1}\cap l_{i_3}\}$, and $\{(x_{i_1}\wedge x_{i_2} \wedge x_{i_3})\mid l_{i_1}\cap l_{i_2}\neq l_{i_1}\cap l_{i_3}\}$.

If $\bk$ is not specified, it will be assumed to be $\bq$.
\end{dfn}

There is a grading in $\ca$ induced by the grading in $E$, and $\ca_1=E_1$. In the following, we will fix the base $\{x_1,\ldots,x_n\}$ to take coordinates.
We will say that two vectors $v_1,v_2\in\ca_1$ are \emph{orthogonal} if $v_1\wedge v_2=0\in\ca$.
\begin{lema}
\label{orto-adm}
Let $\alpha:H\rightarrow\bz^{k-1}$ be a homomorphism, and let $M$ be the matrix whose columns are $\alpha(e_1),\ldots ,\alpha(e_n)$. Then $\alpha$ is admissible if and only if the rows of the matrix $M$ (as elements of $\ca_1$) are orthogonal in $\ca$.
\end{lema}
\begin{proof}
A basis of $\ca_2$ can be given by the generators $x_i\wedge x_j$ such that $l_i$ is the first line of $l_i \cap l_j$, and the rest of the generators can be expressed in terms of these as follows: let $l_i$ be the first line of the point $l_j\cap l_k$, then if $j<k$, we can use $\delta(x_i\wedge x_j\wedge x_k)$ to see that $x_j\wedge x_k=(x_i\wedge x_k)-(x_i\wedge x_j)$. It is not hard to see that the rest of the relations are a consequence of the previous ones.
Now let $a:=(a_1,\ldots ,a_n)$ and $b:=(b_1,\ldots ,b_n)$ be two such rows, $l_j\in p\in \scp$, and let $l_i$ be the first line of $p$. The coefficient of $a\wedge b$ in $x_{i}\wedge x_{j}$ is
\begin{equation}
\label{coefdeter}
\left | \begin{array}{cc} a_{j} & \sum_{l_k\in p}a_{k} \\ b_{j} & \sum_{l_k\in p}b_{k} \end{array} \right|
\end{equation}.
 It is immediate that all such coefficients to be zero is the necessary and sufficient condition for both the admissibility of $\alpha$ and the orthogonality of its rows in $A$.
\end{proof}

\begin{dfn}
For every element $\omega\in E_1$, we consider the complex 
\[
\xymatrix{\bk \ar[r]^{d_\omega} & \ca_1 \ar[r]^{d_\omega}& \ca_2 \ar[r] & 0}
\]
where $d_\omega$ represents the left multiplication by $\omega$. We will denote its cohomology as $H^\bullet(\ca,\omega)$. The \textbf{resonance variety} of $\ca$ is the set $R_1:=\{\omega\in\ca_1\mid H^1(\ca,\omega)\neq 0\}$.
\end{dfn}

\begin{rmrk}
\label{reso-adm}
An element $a\in\ca$ is in $R_1$ if and only if there exists another $b\in\ca\setminus \bk a$, such that $a\wedge b=0$. In that case, the matrix whose rows are $a$ and $b$ defines an admissible map. And vice-versa: the rows of an admissible map are elements of $R_1$. More precisely, maximal admissible classes correspond exactly to irreducible components of the resonance variety \cite{tesis}.
\end{rmrk}

\begin{ntc}
The row vectors of matrix $M$ in Lemma~\ref{orto-adm} will be denoted by $s_1,\ldots,s_n$.
\end{ntc}

We can find a basis $\{r_F\mid F\in\Pi\}$ of the kernel of $Q$ formed by the vectors that generate the kernels of the indecomposable submatrices $M_F$. In particular, we can choose all these vectors to be positive, and to have the property that the sum of their entries is the same for all of them. 
\begin{prop}
The subspace $\langle s_1,\ldots, s_n\rangle_{E_1}$ is exactly $\ker(Q)\cap \ker(U)$.
\end{prop}
\begin{proof}

Equation~\eqref{anulafilas} implies $MQ=0$, since $Q$ is symmetric, $QM^t=0$ and therefore $\langle s_1,\ldots,s_n\rangle\subseteq \ker(Q)$. On the other hand, since $e_1+\cdots +e_n=0$, $\alpha(e_1)+\cdots+\alpha(e_n)=0$, and therefore the sum of the coefficients of each $s_i$ equals zero. Note that $\ker(U)=\{(v_1,\ldots,v_n )\in E_1 \mid \sum_{i=1}^{n}v_i=0\}$ and hence $\langle s_1,\ldots,s_n\rangle\subseteq \ker(Q)\cap \ker(U)$.

For the other inclusion, fix a certain $\bar F\in \Pi$, and define $\tilde r_F:=r_F-r_{\bar F}$ for all $F\in\Pi\setminus\{\bar F\}$. Note that $\{\tilde r_F\}_{F\neq \bar F}$ is a basis of $\ker(Q)\cap \ker(U)$. We will prove that they are pairwise orthogonal (as elements of $\ca$). Take a point $p\in\scp$. If $p\notin \chi_\alpha$, then all the lines in $p$ are in the same component $G$ of $\Pi$. The coefficients of $\tilde r_F\wedge \tilde r_{F'}$ on $x_i\wedge x_j$ are zero for $l_i,l_j \in p$: if $G=\bar F$, the coefficient \eqref{coefdeter} is the determinant of a matrix with two equal rows; if $G=F$ or $G=F'$, one of the rows is zero; and otherwise, both rows are zero.

If $p\in \chi_\alpha$, $\sum_{l\in p}\tilde r_{F,l}$ is the dot product of the row of $J$ corresponding to $p$ by $\tilde r_F$. Since $\tilde r_F$ is in $\ker(Q)\cap \ker(U)$, it is also in $\ker(J^tJ)$. But if $J^tJ\tilde r_F=0$, then $(J\tilde r_F)^t(J\tilde r_F)=\tilde r_F^t J^t J \tilde r_F=\tilde r_F^t\bar 0=0$. Since wall matrices and vectors have real entries $(J\tilde r_F)^t(J\tilde r_F)=0$ implies $J\tilde r_F=0$. Hence, the coefficients of $\tilde r_F\wedge \tilde r_{F'}$ in the generators of $\ca_2$ corresponding to $p$ are again zero.

Therefore, all the $\tilde r_F$'s are pairwise orthogonal. In particular, they are orthogonal to all $s_i$'s, and since $\alpha$ is maximal, the space they span must be the same.
\end{proof}

\begin{thm}
\label{admhaz}
If $\alpha$ is a k-admissible map, then its corresponding admissible subcombinatorics is a combinatorial pencil of no less than k fibers. Furthermore, if $\alpha$ is maximal, then the number of fibers of the pencil is exactly k.
\end{thm}

\begin{proof}

Since there is a basis of $\ker(Q)$ formed by positive vectors, $\ker(Q)\subsetneq \ker(U)$, which means that $\dim(\ker(Q) \cap \ker(U))=\dim(\ker(Q))-1$. This implies that $k=\#\Pi$. 

Equation~\eqref{sumacero} can be expressed as
\begin{equation}
\label{dependenciavectores}
\sum_{F\in\Pi}W_Fv_F=0
,
\end{equation}
where $W_F=\sum_{l\in F}w_l$. We then have a family of $k$ vectors, $\{v_F\}_{F\in \Pi}$ that span $\bz^{k-1}$; that is, all the possible linear combinations satisfied by $\{v_F\}_{F\in \Pi}$ are proportional. By dividing each $v_F$ by a positive integer $z_F$, we may assume that $\sum_{F\in\Pi}v_F=0$. The equations $\{\sum_{l\in p}\alpha(l)=0\mid p\in\chi_\alpha\}$ can be rewritten as $\{\sum_{F\in \Pi}\sum_{l\in F\cap p}w_lz_Fv_F=0\mid p\in\chi_\alpha\}$, which means that at each point $p\in\chi_\alpha$, $\sum_{l\in F\cap p}w_lz_F$ must be constant for all $F\in \Pi$. If we denote $\bar w_l:=z_Fw_l$ for $l\in F\in \Pi$, $(\scl,\scp)$, is a combinatorial pencil with the partition $\Pi$ and the weights $(\bar w_l)_{l\in\scl}$.

\end{proof}
As a direct consequence of the previous Theorem and Remark \ref{reso-adm}, we obtain the following result about the resonance variety:
\begin{thm}
Given a combinatorics $(\scl,\scp)$, there is a bijection between the $k$-dimensional components of its resonance variety and the combinatorial pencils of $k+1$ fibers contained in it.
\end{thm}

\section{Permutations of the admissible classes}
\label{permutaciones}

Let $\Aut^1(H):=\{\phi\in \Aut(H) \mid \phi \wedge \phi(R)=R\}$. From Proposition~\ref{carac-admisible} we obtain the following result.
\begin{cor} 
\label{antihom}
Any $\phi\in Aut^1(H)$ induces a permutation $\sigma_{\phi}$ of the set of $k$-admissible classes by composition. In fact, there is a group antihomomorphism between $Aut^1(H)$ and the group of permutations of $Adm_k(\scl,\scp)$.
\end{cor}
%

These permutations must preserve some structure in the admissible classes. Consider the function $\Upsilon : \calp (Adm(\scl,\scp))\rightarrow \bz$ given by
$$
\Upsilon(S):=\codim(\bigcap_{\alpha\in S} \ker(\alpha)).
$$
Note that if two admissible maps belong to the same admissible class, their kernel must be equal, and hence $\Upsilon$ is well defined. Also note that, if $\alpha$ is a $k$-admissible map, then $\Upsilon(\{\alpha\})=k-1$.

For every $\phi\in Aut^1(H)$, $\sigma_{\phi}$ induces also a permutation $\bar\sigma_\phi$ in $\calp(Adm(\scl,\scp))$. It is straightforward to prove the following.
\begin{lema}
\label{preservaupsilon}
For every $\phi\in Aut^1(H)$, and every $S=\in\calp (Adm(\scl,\scp))$, $\Upsilon(\bar\sigma_\phi(S))=\Upsilon(S)$.
\end{lema}

The previous lemma allows us to calculate the set of the possible $\sigma_\phi$ (which is a subgroup of the permutations of $Adm(\scl,\scp)$, in particular it is the image of the morphism mentioned in Corollary~\ref{antihom}) as follows: consider the natural action of the group of permutations of $Adm(\scl,\scp)$ in $\calp(Adm(\scl,\scp))$. Now, for every couple of positive integers $(i,j)$, consider the subset $P_{i,j}:=\{S\in \calp(Adm(\scl,\scp))\mid \#S=i,\Upsilon(S)=j\}$. Any $\sigma_\phi$ must be in the stabilizer of $P_{i,j}$ for each $(i,j)\in \bz^2$. Therefore, calculating the intersection of all such stabilizers gives us a group that contains $\{\sigma_\phi\mid \phi\in Aut^1(H)\}$ as a subgroup.

But in most cases it is enough to use a particular version of the previous method, by considering only the concept of triangle, which we define below.

\begin{dfn}

Let $\alpha_1$, $\alpha_2$ and $\alpha_3$ be admissible maps. We will say that they form \textbf{a triangle} (of admissible classes) if 

\begin{equation}
\Upsilon(\{\alpha_1, \alpha_2,\alpha_3\})=\sum_{i=1}^3\Upsilon(\{\alpha_i\})-1
\end{equation}

If the admissible classes are clear from the context, we can also talk about triangles of admissible subcombinatorics.
\end{dfn}

\begin{ejm}
Let $p_1$, $p_2$ and $p_3$ be three points of multiplicities $m_1$, $m_2$ and $m_3$ greater than two. They are admissible subcombinatorics as in Example~\ref{punto}. These admissible classes form a triangle if and only if $p_i \cap p_j \neq \varnothing$ for $i,j =1,2,3$ but $p_1 \cap p_2 \cap p_3 = \varnothing$. This fact can be easily checked case by case (see \cite{nuestro-ryb} for details).
\end{ejm}

A direct consequence of Lemma~\ref{preservaupsilon} is the following.
\begin{lema}
\label{preservatriangulos}
If $\alpha_1$, $\alpha_2$ and $\alpha_3$ are $k_1$, $k_2$ and $k_3$-admissible maps respectively that form a triangle, their images under $\sigma_\phi$ are $k_1,k_2$ and $k_3$ admissible maps that form a triangle too. Thus $\sigma_\phi$ maps triangles to triangles.
\end{lema}

This fact can be used to study the possible permutations induced by $\phi$ and, in certain cases, that is enough to calculate all the possible automorphisms of $H$ that fix $R$.

\section{An example}
\label{ejemplo}

For some combinatorics $\Aut^1(H)$ is as small as it gets, that is $\{\pm Id\}\times \Aut(\scl,\scp)$ where $\Aut(\scl,\scp)$ is the group of automorphisms of the combinatorics. Such combinatorics are called \textbf{homologically rigid}.

The existence of two real line arrangements with the same combinatorial type, but different topology of the embedding was shown in \cite{nuestro-oncerectas}. These arrangements have 11 lines, and are conjugated in $\bq (\sqrt{5})$. They can be constructed after two other arrangements (also conjugated in $\bq (\sqrt{5})$), whose non-generic braid monodromies were shown to be non-equivalent. Here we will show that the combinatorics of these arrangements of ten lines is homologically rigid  by studying the possible permutations of the admissible classes induced by $\Aut^1(H)$. This combinatorics is the result of adding one line to the Falk-Sturmfels combinatorics described in~\cite{cohen-suciu-braidmonodromy}. This example is intended to show how to use Lemmas \ref{preservaupsilon} and \ref{preservatriangulos} to study the homological rigidity of a given combinatorics. 

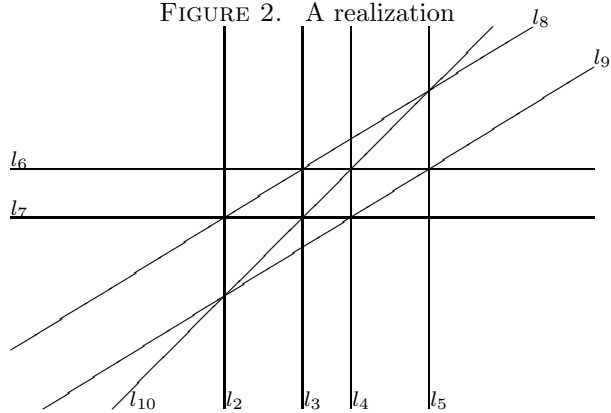
\begin{figure}[h]
\caption{\label{diezrectas} A realization}
\setlength{\unitlength}{0.00087489in}
\begingroup\makeatletter\ifx\SetFigFont\undefined%
\gdef\SetFigFont#1#2#3#4#5{%
  \reset@font\fontsize{#1}{#2pt}%
  \fontfamily{#3}\fontseries{#4}\fontshape{#5}%
  \selectfont}%
\fi\endgroup%
{\renewcommand{\dashlinestretch}{30}
\begin{picture}(3534,2339)(0,-10)
\drawline(12,1162)(3522,1162)
\drawline(12,1454)(3522,1454)
\drawline(1767,2312)(1767,12)
\drawline(2059,2312)(2059,12)
\drawline(2914,2312)(623,12)
\drawline(2529,2312)(2529,12)
\drawline(3151,2312)(12,363)
\drawline(1299,2312)(1299,12)
\drawline(3522,2070)(206,12)
\put(1309,12){\makebox(0,0)[lb]{{\SetFigFont{8}{10}{\rmdefault}{\mddefault}{\updefault}$l_2$}}}
\put(1777,12){\makebox(0,0)[lb]{{\SetFigFont{8}{10}{\rmdefault}{\mddefault}{\updefault}$l_3$}}}
\put(2069,12){\makebox(0,0)[lb]{{\SetFigFont{8}{10}{\rmdefault}{\mddefault}{\updefault}$l_4$}}}
\put(2539,12){\makebox(0,0)[lb]{{\SetFigFont{8}{10}{\rmdefault}{\mddefault}{\updefault}$l_5$}}}
\put(12,1454){\makebox(0,0)[lb]{{\SetFigFont{8}{10}{\rmdefault}{\mddefault}{\updefault}$l_6$}}}
\put(12,1162){\makebox(0,0)[lb]{{\SetFigFont{8}{10}{\rmdefault}{\mddefault}{\updefault}$l_7$}}}
\put(730,12){\makebox(0,0)[lb]{{\SetFigFont{8}{10}{\rmdefault}{\mddefault}{\updefault}$l_{10}$}}}
\put(3151,2312){\makebox(0,0)[lb]{{\SetFigFont{8}{10}{\rmdefault}{\mddefault}{\updefault}$l_8$}}}
\put(3522,2070){\makebox(0,0)[lb]{{\SetFigFont{8}{10}{\rmdefault}{\mddefault}{\updefault}$l_9$}}}
\end{picture}
}

\end{figure}

The combinatorics that we will study has ten lines $\{l_1,\ldots ,l_{10}\}$, ten triple points: $\{l_1,l_6,l_7\}$, $\{l_1,l_8,l_9\}$, $\{l_2,l_9,l_{10}\}$, $\{l_2,l_7,l_8\}$, $\{l_3,l_6,l_8\}$, $\{l_3,l_7,l_{10}\}$, $\{l_4,l_6,l_{10}\}$, $\{l_4,l_7,l_9\}$, $\{l_5,l_8,l_{10}\}$, and $\{l_5,l_6,l_9\}$, and a quintuple point $\{l_1,l_2,l_3,l_4,l_5\}$. The remaining are double points. A real realization, where $l_1$ is the line at infinity, can be seen in Figure~\ref{diezrectas}.

 We could calculate the admissible classes by solving the quadratic equation system that the coefficients of any admissible map should satisfy, but it is much faster to use Theorem~\ref{admhaz}, which allows us to calculate the possible combinatorial pencils just by solving systems of linear equations. The result is that the only non-point-type admissible classes are the following ten Ceva type classes: 
 \begin{itemize}
\item $\{l_1,l_4,l_5,l_9,l_6,l_7\}$
\item  $\{l_2,l_4,l_5,l_6,l_9,l_{10}\}$
\item  $\{l_1,l_2,l_4,l_7,l_9,l_8\}$
\item  $\{l_1,l_2,l_3,l_8,l_6,l_7\}$\item  $\{l_2,l_3,l_5,l_{10},l_8,l_7\}$\item  $\{l_1,l_3,l_4,l_{10},l_6,l_7\}$\item  $\{l_3,l_4,l_5,l_{10},l_8,l_6\}$\item $\{l_1,l_2,l_5,l_{10},l_8,l_9\}$\item  $\{l_1,l_3,l_5,l_6,l_9,l_8\}$ \item $\{l_2,l_3,l_4,l_7,l_9,l_{10}\}$
 \end{itemize}

There is only one 5-admissible subcombinatorics, which must be preserved by the permutation induced by any $\phi\in \Aut^1(H)$. For each of the 20 3-admissible classes, we can count to how many triangles of maximal 3-admissible maps it belongs. The result is that each point-type combinatorics belongs to 15 such triangles, while each Ceva type belongs to 9. Hence $\sigma_{\phi}$ must induce a permutation of the triple points that preserves triangles. These computations were done in a few seconds in a computer using GAP \cite{gap}.

Now consider four points $p_1,p_2,p_3,p_4$ of multiplicity greater than two such that $p_i,p_j,p_4$ form a triangle for all $i,j\in \{1,2,3\}$, $i\neq j$. Then $p_1,p_2$ and $p_3$ are aligned if and only if they do not form a triangle. In our combinatorics, for any three aligned points of multiplicity greater than two, there is a fourth one that forms a triangle with any two of them. Since $\sigma_\phi$ must preserve triangles, it must also preserve lines considered as sets of multiple points. Up to composition with the automorphism of $H$ induced by $\sigma_\phi^{-1}$ (seen as a permutation of the elements of the canonical generating system of $H$), we may assume that $\sigma_\phi$ is the identity.

Now we can consider the basis of $H$ given by $\bar{e}_1,\ldots ,\bar{e}_9$ (being $\bar{e}_i$ the class of $e_i$ modulo $(1,\ldots, 1)$) and the matrix $A=(a_{i,j})$ related to $\phi$ in this basis. Using the fact that $\sigma_\phi$ is the identity, we can deduce that for any point $p$ of multiplicity greater than two, $a_{{i},k}=a_{{j},k}$ for all $l_{i},l_{j} \in p$ and $l_k\notin p$. If $l_{10}\in p$ then these entries are actually 0. These conditions to all the multiple points, forces $A$ to be diagonal. Since $A$ must be an integer matrix with determinant equal to $\pm 1$, all the entries in the diagonal must be $\pm 1$. Now given a multiple point $p$ such that $l_{10}\notin p$, we have that the submatrix of $A$ obtained by selecting the rows and columns corresponding to the lines in $p$ must have columns that add up to a multiple of $(1,\ldots ,1)$, so it means that $a_{i,i}=a_{j,j}$ for all $l_i,l_j\in p$. If we use these conditions in all the multiple points, we obtain that $A$ must be $\pm Id$. 

This method also works with the combinatorics of McLane and Rybnikov \cite{rybnikov}, and the one of eleven lines studied in \cite{nuestro-oncerectas}. This same kind of arguments were used by Falk in \cite{falk-arrco} to show that certain family of combinatorics (which he called strongly connected) are homologically rigid. In fact this method could be seen as a generalization of his.

\section{Duality between $(H \wedge H)/R$ and the Orlik-Solomon algebra}
\label{dualidad}

The method in the previous section can allow us to calculate $\Aut^1(H)$. Here we will see that this group coincides with the group of automorphisms of the Orlik-Solomon algebra.

Let $A=(a_{i,j})$ be the matrix that represents an automorphism of $H$ in the basis $\{e_1,\ldots ,e_n\}$. As a sublattice of $H\wedge H$, $R$ is generated by $\{e_i\wedge \sum_{l_j\in p} e_j \mid l_i \in p , p\in  \scp\}$. Since $\sum_{l_i\in p} (e_i\wedge \sum_{l_j\in p}e_j)=0$ when $\# p >2$, we can eliminate the first one and then use these relations to give a basis of the quotient. In particular, for each point $p=\{l_{i_1},\ldots,l_{i_{\# p}}\}$, we can express the generators of the form $e_{i_1}\wedge e_{i_j}$ as $\sum_{k=j+1}^{\# p}e_{i_j}\wedge e_{i_k} - \sum_{k=1}^{j-1} e_{i_k}\wedge e_{i_j}$. Therefore a basis of the quotient is given by the generators of the form $e_i\wedge e_j$ where $i<j$ and both $i$ and $j$ are not the first line in $l_i \cap l_j$.

The image of $e_{i_1}\wedge e_{i_2}$ under the automorphism induced by $A$ is $$\sum_{1\leq j_1<j_2\leq n} \left | \begin{matrix}
a_{j_1,i_1} & a_{j_1,i_2} \\ 
a_{j_2,i_1} & a_{j_2,i_2}
\end{matrix} \right |(e_{j_1}\wedge e_{j_2}),
 $$so the image of a generator of $R$ of the form $e_i\wedge \sum_{l_k\in p} e_k$ is $$\sum_{1\leq j_1<j_2\leq n} \left | \begin{matrix}
a_{j_1,i} & \sum_{l_k\in p} a_{j_1,k} \\ 
a_{j_2,i} & \sum_{l_k\in p} a_{j_2,k}
\end{matrix} \right |(e_{j_1}\wedge e_{j_2}).
 $$
If we project this to the quotient, its coefficient in a generator of the form $e_c\wedge e_d$ (if the first line that goes through the intersection point of $l_c$ and $l_d$ is $l_b$) is
\begin{equation}
\label{resulr2}
\left | \begin{matrix}
a_{c,i} & \sum_{l_k\in p} a_{c,k} \\ 
a_{d,i} & \sum_{l_k\in p} a_{d,k}
\end{matrix} \right |
-
\left | \begin{matrix}
a_{b,i} & \sum_{l_k\in p} a_{b,k} \\ 
a_{d,i} & \sum_{l_k\in p} a_{d,k}
\end{matrix} \right |
+
\left | \begin{matrix}
a_{b,i} & \sum_{l_k\in p} a_{b,k} \\ 
a_{c,i} & \sum_{l_k\in p} a_{c,k}
\end{matrix} \right |
=
\left | \begin{matrix}
a_{b,i} & \sum_{l_k\in p} a_{b,k}  & 1\\ 
a_{c,i} & \sum_{l_k\in p} a_{c,k} & 1\\
a_{d,i} & \sum_{l_k\in p} a_{d,k} & 1
\end{matrix} \right | .
\end{equation}

The extra condition that $\phi_A\in\Aut^1(H)$ is equivalent to asking \eqref{resulr2} to vanish for every point $p$, every line $l_i \in p$ and every three concurrent lines $l_b,l_c,l_d$ such that $l_b$ is the first line in their intersection point.

Now let's look at the second level of the Orlik-Solomon algebra $\ca^2$. The relations we have are of the form $(x_b\wedge x_c)-(x_b\wedge x_d)+(x_c\wedge x_d)$ for every three concurrent lines $l_b,l_c,l_d$. Now let's suppose that there are three concurrent lines such that the first line that goes through their intersection point is $l_a$; we can express the relation $(x_b\wedge x_c)-(x_b\wedge x_d)+(x_c\wedge x_d)$ as $((x_a\wedge x_b)-(x_a\wedge x_c)+(x_b\wedge x_c))$\(-\)$((x_a\wedge x_b)-(x_a\wedge x_d)+(x_b\wedge x_d))$\(+\)$((x_a\wedge x_c)-(x_a\wedge x_d)+(x_c\wedge x_d))$. In particular we only need the relations where the first line in the intersection point appears. These relations allow us to express any $x_b\wedge x_c$ as $(x_a\wedge x_c)-(x_a\wedge x_b)$ (where again $b<c$ and $l_a$ is the first line that goes through $l_b \cap l_c$. Hence a basis of $\ca^2$ is given by the $x_i\wedge x_j$ such that $l_i$ is the first line in $l_i \cap l_j$. Let $l_b,l_c,l_d$ be three concurrent lines; the image of $(x_b\wedge x_c)-(x_b\wedge x_d)+(x_c\wedge x_d)$ by $\phi_A\wedge \phi_A$ is 
\begin{equation}
\label{imagenorliksolomon}
\begin{array}{l}
\sum_{1\leq i < j \leq n} \left( \left | \begin{matrix}
a_{i,b} & a_{i,c} \\ 
a_{j,b} & a_{j,c} 
\end{matrix} \right |  -
\left | \begin{matrix}
a_{i,b} & a_{i,d} \\ 
a_{j,b} & a_{j,d} 
\end{matrix} \right |
+  \left | \begin{matrix}
a_{i,c} & a_{i,d} \\ 
a_{j,c} & a_{j,d} 
\end{matrix} \right |
\right)
 (x_i\wedge x_j)
= \\
=\sum_{1\leq i < j \leq n} \left | \begin{matrix}
a_{i,b} & a_{i,c} & a_{i,d} \\ 
a_{j,b} & a_{j,c} & a_{j,d} \\
1 & 1 & 1
\end{matrix} \right |
(x_i\wedge x_j).
\end{array}
\end{equation}

Now let's take a point $p$ whose first line is $l_i$ and $l_j,l_k\in p$; by the previous relations, $x_j\wedge x_k$ maps by the projection to $x_i\wedge x_k-x_i\wedge x_j$. This means that if we want to calculate the coefficient of the projection of \eqref{imagenorliksolomon} in $x_i\wedge x_j$, we have to add or substract adequately its coefficients in all $x_j\wedge x_k$ such that $l_k \in p$. The result is that the coefficient of \eqref{imagenorliksolomon} in $x_i\wedge x_j$ is 

\begin{equation}
\label{resulorliksolomon}
\left | \begin{matrix}
\sum_{k\in p} a_{k,a} & \sum_{k\in p}a_{k,b} & \sum_{k\in p}a_{k,c} \\ 
a_{j,a} & a_{j,b} & a_{j,c} \\
1 & 1 & 1
\end{matrix} \right |
.
\end{equation}

In order for $A$ to induce an automorphism of the Orlik-Solomon algebra, \eqref{resulorliksolomon} must hold for all three concurrent lines $l_b,l_c,l_d$, each point $p$, and each line $l_j\in p$. Comparing \eqref{resulr2} and \eqref{resulorliksolomon}, we have the following result.
\begin{prop}
The matrix $A$ induces an automorphism of the Orlik-Solomon algebra if and only if its transposed $A^t$ induces an automorphism of $H\wedge H$ that preserves $R$.
\end{prop}

\bibliography{biblio} 

\end{document}